\documentstyle[amssymb]{article}
\input BoxedEPS
\SetepsfEPSFSpecial
\begin{document}
\title{The Santal${\mbox{\'o}}$-regions of a convex body\thanks{ the paper
was written
while both authors stayed at MSRI}}
\author{Mathieu Meyer and Elisabeth Werner\thanks{ supported by a
grant from the National Science Foundation.
\newline
\newline
MSC classification   52
}}

\date{}
\maketitle

\begin{abstract}
Motivated by the Blaschke-Santal${\mbox{\'o}}$ inequality, we define for a
convex
body K in ${\bf R}^n$ and for $t \in {\bf R}$ the
Santal${\mbox{\'o}}$-regions S(K,t) of K. We investigate properties of these
sets and  relate them to a concept of Affine Differential Geometry,
the affine surface area of K.
\end{abstract}

\vskip 10mm
Let K be a convex body in ${\bf R}^n$. For $x\in$ int(K), the
interior of K, let $K^x$ be the polar body of K with respect
to x. It is well known that there exists a unique $x_0 \in$
int(K) such that the product of the volumes $|K||K^{x_0}|$ is
minimal (see for instance [Sch]). This unique $x_0$ is called the
Santal${\mbox{\'o}}$-point of K.
\par
Moreover the Blaschke-Santal${\mbox{\'o}}$ inequality says that
$|K||K^{x_0}| \leq v_n^2$ (where $v_n$ denotes the volume of the
n-dimensional Euclidean unit ball $B(0,1)$) with equality if and
only if K is an ellipsoid.
\vskip 5mm
For $t \in {\bf R}$ we consider here the sets
$$S(K,t)=\{x \in K : \frac{|K||K^{x}|}{v_n^2} \leq t\}.$$
Following E. Lutwak, we call $S(K,t)$ a Santal${\mbox{\'o}}$-region of K.
\newline
Observe that it follows from the Blaschke-Santal${\mbox{\'o}}$ inequality that
the Santal${\mbox{\'o}}$-point $x_0 \in S(K,1)$ and that $S(K,1) = \{x_0\}$
if and only if K is an ellipsoid. Thus $S(K,t)$ has non-empty interior
for some $t < 1$ if and only if K is not an ellipsoid.
\vskip 5mm
In the first part of this paper we show some properties of
$S(K,t)$ and give estimates
on the ``size" of $S(K,t)$. This question was asked by E. Lutwak.
\vskip 5mm
In the second part we show how $S(K,t)$ is related to the affine surface
area of K.
\par
The affine surface area $as(K)$ is originally a notion of
Differential Geometry.
For a convex body $K$ in ${\bf R}^n$ with sufficiently smooth boundary
$\partial K$ it
is
defined as
$$as(K)=\int_{\partial K} {\kappa} (x)^{\frac{1}{n+1}} d \mu (x),$$
where ${\kappa} (x)$ is the Gaussian curvature in $x \in \partial K$ and
$\mu$ is the surface measure on $\partial K$.
The affine surface area is invariant under
affine transformations with determinant 1. It arises naturally in
questions concerning the
approximation of convex bodies by polytopes (see [G]) and in a
priori estimates of PDE's ([Lu-O]).
\par
It has been one of the aims of Convexity Theory to extend the notions of
Differential Geometry and for instance of affine surface area
to arbitrary convex bodies without any
smoothness assumptions on the boundary.
\par
Within the last few years  four different extensions
have been
given (due to Leichtweiss [L1],  Lutwak [Lu], Sch\"{u}tt-Werner [S-W] and
Werner[W]) and it was shown that they all coincide  ([S1], [D-H]).
\par
We give here another such extension arising again from a completely different
context. It will also follow that this new extension coincides with the others.
\vskip 5mm
The authors wish to thank MSRI for the hospitality and the organizers of
the special semester in Convex Geometry  and Geometric Functional
Analysis at MSRI for inviting them. It was during the stay there that the
paper was written.

\vskip 10mm
Unless stated otherwise we will always assume that a convex body K
in ${\bf R}^n$ has its Santal${\mbox{\'o}}$-point
at the origin. Then 0 is the center of mass of the polar body $K^0$
which may be written as
$$\int_{K^0}<x,y> dy = 0
\hspace{.1in}\mbox{for every}\hspace{.1in} x \in {\bf R}^n.$$
By $|K|$ we denote the n-dimensional volume of K. $h_K$ is the
support function of K. If K is   centrally symmetric, $||.||_K$ is the norm on
${\bf R}^n$ that  has K as its unit ball. By $||.||$ we denote the standard
Euclidean  norm on ${\bf R}^n$, $<^.,^.>$ is the usual inner product on
${\bf R}^n$. $B(a,r)$ is the n-dimensional
Euclidean ball with radius r
centered at a. For $x\in K$, $K^x = (K-x)^0=\{y\in {\bf R}^n: <y,z-x> \leq 1
\hspace{.1in}\mbox{for all z} \in K \}$ is the polar body of K with respect
to x;
$K^0$
denotes the polar body with respect to the Santal${\mbox{\'o}}$-point. Moreover
for
$u \in S^{n-1}$ we will denote by $\phi_K^u(y)$ or in short by $\phi(y)$
the $(n-1)$ - dimensional volume of the sections of K orthogonal to u, that is
$$ \phi(y)=\phi_K^u(y) =  |\{z \in K: <z,u>=y\}|.$$
\newpage
I. PROPERTIES OF THE SANTAL${\mbox{\'O}}$-REGIONS
\vskip 5mm
We start by listing some of the properties of $S(K,t)$.
Recall that for $\delta > 0$, $\delta$ small enough, $K_{\delta}$
is said to be a (convex) floating body of $K$, if it is the intersection of all
halfspaces whose defining hyperplanes cut off a set of volume
$\delta$ of K ([S-W]). More precisely, for $ u \in S^{n-1}$ and for
$ 0 < \delta$ let $a_{\delta}^u$ be defined by
$$|\{x \in K: <x,u> \geq a_{\delta}^u\}| = \delta.$$ Then $K_{\delta}
= \cap_{u \in S^{n-1}}\{x\in K: <x,u> \leq a_{\delta}^u\}$.

\vskip 5mm
In the following proposition we consider only those $t \in {\bf R}$ for which
\newline
$S(K,t) \neq \emptyset$.
\vskip 3mm
\medskip\noindent {\bf Proposition 1}

\begin{it}
Let K be a convex body in ${\bf R}^n$. Then
\vskip 2mm
(i) $S(K,t)$ is strictly convex for all t.
\vskip 2mm
(ii) $S(A(K),t) = A(S(K,t))$ for all regular affine transformations A,
\newline
for all t.
\vskip 2mm
(iii) The boundary of $S(K,t)$ is $C^{\infty}$ for all t.
\vskip 2mm
(iv) $t\longmapsto S(K,t)$ is increasing and concave, that is for all t, s
and for all $\alpha \in {\bf R}$, $0 \leq \alpha \leq 1,$
$$S(K,\alpha t + (1-\alpha ) s) \supset \alpha S(K,t)
+ (1- \alpha )S(K,s).$$
\vskip 2mm
(v) For all $0 < \delta < \frac{1}{2}$,   $K_{\delta
|K|}$  is contained in $S(K,\frac{1}{4 \delta (1-\delta)})$.
\vskip 5mm
\end{it}
\medskip\noindent {\bf Proof}

We frequently use the following well known formula.
For all $x \in \mbox{int}(K)$
\begin{equation}
 |K^x|= \frac{1}{n} \int_{S^{n-1}} \frac{d\sigma(u)}
{(h_K(u) - <u,x>)^n} ,
\end{equation}
where $\sigma$ is the spherical Lebesgue measure.
\newline
Indeed, $|K^x|= \frac{1}{n} \int_{S^{n-1}} \frac{d\sigma(u)}
{(\tilde{h}_K(u) )^n}$, where $\tilde{h}_K$ is the support function of $K$
centered at x. Now observe that $\tilde{h}_K(u) =
h_K(u) - <u,x>$; thus (1) follows.
\vskip 5mm
(i) Observe that for all $u \in S^{n-1}$ the function
$$x \longmapsto \frac{1}{(h_K(u) -
<x,u>)^n} $$
is convex on int(K). (1) then implies that
$$x\longmapsto |K^x|=  \frac{1}{n} \int_{S^{n-1}} \frac{d\sigma(u)}
{(h_K(u) - <u,x>)^n}$$
is convex. In fact the function
$$x\longmapsto |K^x|=  \frac{1}{n} \int_{S^{n-1}} \frac{d\sigma(u)}
{(h_K(u) - <u,x>)^n}$$
is strictly convex, as for $x,y \in \mbox{int}(K), x\neq y$,
$$\sigma (\{u \in S^{n-1}: <u,x> = <u,y> \}) = 0.$$
Therefore (i) follows.
\vskip 5mm
(ii) Let A be a one-to-one affine transformation. We can write $A=L+a$,
where L is a one-to-one linear transformation and $a$ is a vector in
${\bf R}^n$.  Then
$$(A(K))^y = \{x \in {\bf R}^n: <x, Az - y> \leq 1\hspace{.1in} \mbox{for all
z} \in K\}$$
$$=\{x \in    {\bf R}^n: <L^*x, z -A^{-1} y> \leq 1\hspace{.1in} \mbox{for all
z} \in K\}$$
$$=(L^{*})^{-1}(K^{A^{-1}(y)}).$$
Hence
$$S(AK,t) = \{y \in AK: \frac{|AK||AK^y|}{v_n^2} \leq t\}$$
$$= \{y \in AK: \frac{|K||K^{A^{-1}(y)}|}{v_n^2} \leq t\}$$
$$=A(S(K,t)).$$
\vskip 5mm
(iii) Let K, n and t be fixed. By (1), $\partial
S(K,t)=\{x\in K: F(x)=\frac{n t v_n^2}{|K|}\}$, where
$F(x)=\int_{S^{n-1}}
\frac{d\sigma(u)} {(h_K(u) - <u,x>)^n}$. $F$ is continuous on int(K) with
continuous partial derivatives of all orders, has a unique  minimum
at the Santal${\mbox{\'o}}$-point $x_0$ and is convex (see (i)). Therefore
(iii) follows
from the implicit function theorem.
\vskip 5mm
(iv) is obvious from the proof of (i).
\vskip 5mm
(v) Let $\delta \in (0,\frac{1}{2})$, let H be a hyperplane that has
non-empty intersection with K and is such that
$$|K \cap H^+| =\delta |K|,$$
where $H^+$ is one of the two halfspaces determined by H.
By definition the convex floating body $K_{\delta|K|}$ is the
intersection of all the halfspaces $H^-$ determined by all such
hyperplanes H.
\par
On the other hand by [Me-P] there exists $x\in H \cap \mbox{int}(K)$ such
that
$$\frac{|K||K^{x}|}{v_n^2} \leq \frac{1}{4\delta(1-\delta)}.$$
This means that  $x \in S(K,\frac{1}{4\delta(1-\delta)})$.
Consequently
$$K_{\delta|K|} \subseteq S(K,\frac{1}{4\delta(1-\delta)}).$$
\vskip 5mm
\medskip\noindent {\bf Remark 2}
\par
(i) Proposition 1,(iv) says that
$K_{\delta|K|} \subseteq S(K,\frac{1}{4\delta(1-\delta})$.
We will show
(see Propositon 14 ) that in the case of a convex body with
sufficiently smooth boundary and positive Gaussian curvature everywhere
a converse inclusion holds for
$\delta$ ``small".
\vskip 3mm
(ii) Note also that for $K=B(0,1)$,

$$(B(0,1))_{\delta v_n} \sim S(B(0,1),\frac{v_{n-1}}{\delta(n+1)v_n}),$$
for $\delta$ sufficiently small. More precisely, for
$\delta \leq \frac{2^{n+1}v_{n-1}}{\sqrt e(n+1)v_n n^{\frac{n+1}{2}}}$
$$ S(B(0,1),\frac{v_{n-1}}{\sqrt e\delta(n+1)v_n}) \subseteq
(B(0,1))_{\delta v_n} \subseteq
S(B(0,1),\frac{(\sqrt e)^{\frac{n+1}{n-1}}v_{n-1}}{\delta(n+1)v_n}).$$

This follows from the forthcoming Corollary 5 and from the fact
that the volume of a cap of the Euclidean unit ball of height $\Delta$
can be estimated from above by  $\frac{v_{n-1}}{n+1} (2
\Delta)^{\frac{n+1}{2}}$
and from below by $\frac{v_{n-1}}{n+1} (2 \Delta)^{\frac{n+1}{2}}
(1-\frac{\Delta}{2})^{\frac{n-1}{2}}.$
\vskip 3mm
For $\delta$  ``close" to $\frac{1}{2}$,
$$ (B(0,1))_{\delta v_n} \sim
S(B(0,1),\frac{1}{4\delta(1-\delta)}).$$
More precisely, let $\varepsilon \leq \frac{1}{\sqrt n}$ and
$\frac{1}{2} > \delta \geq \frac{1}{2}
-\varepsilon \frac{v_{n-1}}{\sqrt e\hspace{0.1in} v_n}$.
Then
$$ (B(0,1))_{\delta v_n} \subseteq
S(B(0,1),\frac{1}{4\delta(1-\delta)}) \subseteq
\frac{8\hspace{0.1in} v_{n-1}}{\sqrt {n+1}\hspace{0.1in}
v_n}(B(0,1))_{\delta v_n}.$$
\vskip 5mm
The following lemmas will enable us to compute $|(B(0,1))^x|$.
They are also needed for Part II.

\vskip 5mm
\medskip\noindent {\bf Lemma 3}
\par
\begin{it}
Let $ x \in \mbox{int}(K)$. Then
$$|K^x| = \int_{K^0} \frac{dy}
{(1-<x,y>)^{n+1}}.$$
\end{it}
\vskip 3mm
\medskip\noindent {\bf {Proof}}
\par
By (1)
$$|K^x|=\frac{1}{n} \int_{S^{n-1}} \frac{d\sigma(u)}
{(h_K(u))^n(1 - <\frac{u}{h_K(u)},x>)^n}$$
$$=\int_{S^{n-1}} \int_{0}^{\frac{1}{h_K(u)}}\frac{r^{n-1}}{(1 <x,ru>)^{n+1}}
drd\sigma(u)=\int_{K^0} \frac{dy}
{(1-<x,y>)^{n+1}}.$$

\vskip 6mm
\medskip\noindent {\bf Remark}
\newline
We will use Lemma 3 mostly in the following form:
\newline
let $ u \in S^{n-1}$ and $\lambda \in {\bf R}$ such that $x=\lambda u \in
\mbox{int}(K)$. Then  Lemma 3 says that
\begin{equation}
|K^x| = \int_{-h_{K^0}(-u)}^{h_{K^0}(u)} \frac{\phi_{K^0}^u(t)
dt} {(1-\lambda t)^{n+1}},
\end{equation}
where $\phi_{K^0}^u(t) =  |\{z \in K^0: <z,u>=t\}|$.
\vskip 6mm
\medskip\noindent {\bf Lemma 4}
\par
\begin{it}
(i) Let $0 \leq \alpha < 1$. Then
$$\int_{-1}^1 \frac{(1-x^2)^{\frac{n-1}{2}} dx}{(1-\alpha x)^{n+1}}
= \frac {2^n (\Gamma(\frac{n+1}{2}))^2}{(1-\alpha^2)^{\frac{n+1}{2}}
n!}$$
\par
(ii) For $\alpha \in (0,1)$ let
$$I(\alpha) =
(\int_{0}^1 \frac{(1-x^2)^{\frac{n-1}{2}} dx}{(1-\alpha x)^{n+1}})
\hspace{.1in}(\frac
{\alpha^{\frac{n+1}{2}}(1-\alpha)^{\frac{n+1}{2}}
n!}{2^{\frac{n-1}{2}}
(\Gamma(\frac{n+1}{2}))^2}).$$ Then
$$I(\alpha ) \leq 1 \hspace{.1in} \mbox{and}\hspace{.1in}
\lim_{\alpha \rightarrow 1}I(\alpha ) = 1.$$
\par
(iii) Let $a,b >0, n \in {\bf N}$ and $\lambda a<1$. Then
$$\int_{-b}^a \frac{(1-\frac{y}{a})^{n-1}}{(1-\lambda y)^{n+1}}dy
=\frac{(a+b)^n}{na^{n-1}(1-\lambda a)(1+\lambda b)^n}.$$

\end{it}

\medskip\noindent {\bf {Proof}}
\par
(i) We put
$x=\frac{1-\frac{1-\alpha}{1+\alpha}u}{1+\frac{1 \alpha}{1+\alpha}u}$.
This gives (i).
\vskip 5mm
(ii) Put $x=1-w\frac{1-\alpha}{\alpha}$. Then
$$I(\alpha) =\frac{n!}
{2^{\frac{n-1}{2}}
(\Gamma(\frac{n+1}{2}))^2}
\int_0^{\frac{\alpha}{1-\alpha}}\frac{w^{\frac{n-1}{2}}
(2-\frac{1-\alpha}{\alpha}w)^{\frac{n-1}{2}}dw}{(1+w)^{n+1}}.$$
The upper estimate for (ii) follows immediately from this
last expression.
\newline
And by the Monotone Convergence Theorem this last expression tends to
$$\frac{n!}
{
(\Gamma(\frac{n+1}{2}))^2}\int_0^{\infty}\frac{w^{\frac{n-1}{2}}
dw}{(1+w)^{n+1}},$$
which is equal to 1.
\vskip 5mm
(iii) Note that
$$\int \frac{(1-\frac{y}{a})^{n-1}}{(1-\lambda y)^{n+1}}dy =
\frac{(1-\frac{y}{a})^n}{n(\lambda-\frac{1}{a})(1-\lambda y)^n}.$$

This immediately implies (iii).
\vskip 7mm
\medskip\noindent {\bf Corollary 5}
\par
\begin{it}
Let $B(0,r)$ be the n-dimensional Euclidean ball with radius r centered at 0.
For $u \in S^{n-1}$ let $x=\lambda u$, $0 \leq \lambda <r$.
Then
$$|(B(0,r))^x| =
\frac{v_n}{r^n(1-(\frac{\lambda}{r})^2)^{\frac{n+1}{2}}}$$
\end{it}
\par
\medskip\noindent {\bf {Proof}}
\par
The proof follows from (2) and Lemma 4 (i).
\vskip 7mm
Next we estimate the ``size" of $S(K,t)$ in terms of ellipsoids.
Recall that for a convex body K the Binet ellipsoid $E(K)$
is defined by (see[Mi-P])
$$||u||_{E(K)}^2 = \frac{1}{|K|}\int_K <x,u>^2dx , \hspace{.3in}  \mbox{for
all u}
\in{\bf R}^n.$$ We first treat the case when K is a symmetric convex body.
\newpage
\medskip\noindent {\bf Theorem 6}
\par
\begin{it}
Let K be a symmetric convex body in ${\bf R}^n$. For all $t \in {\bf R}$

$$d_{n}(t) \hspace{.1in}E(K^0)
\subseteq S(K,t) \subseteq c_{n}(t)\hspace{.1in}E(K^0),$$

where $$d_{n}(t) =\frac{1}{\sqrt3
\hspace{.1in}n}\hspace{.1in}(1-\frac{|K||K^0|}{tv_n^2})^{\frac{1}{2}}$$

and

$$c_{n}(t) =\mbox{min}\{(\frac{2}{(n+1)(n+2)})^{\frac{1}{2}}
(\frac{tv_n^2}{|K||K^0|})^{\frac{1}{2}}
\hspace{.1in}(1-\frac{|K||K^0|}{tv_n^2})^{\frac{1}{2}},\hspace{.1in}
\sqrt2
(1-(\frac{|K||K^0|}{tv_n^2})^{\frac{1}{n}})^{\frac{1}{2}} \}$$

\end{it}
\vskip 5mm
\medskip\noindent{\bf {Remarks}}
\par
(i) Especially for any ellipsoid $E$, $S(E,1) = \{0\}$.
\vskip 3mm
(ii) If $t \rightarrow \frac{|K||K^0|}{v_n^2}$, then  $S(K,t) \rightarrow
\{0\}$.
\vskip 3mm
(iii) The second expression in $c_{n}(t)$ gives a better estimate from
above than the first
iff
$\frac{|K||K^0|}{tv_n^2}$ is of a smaller order of magnitude than $(n
\mbox{log}n)^{-1}$.
\vskip 3mm
(iv) Recall that for two isomorphic Banach spaces E and F
the Banach-Mazur distance $d(E,F)$ is defined by
$$d(E,F) = \mbox{inf}\{||T||||T^{-1}||: \mbox{T is an isomorphism
from E onto F}\}.$$
For symmetric convex bodies $K$, $L$ in ${\bf R}^n$ denote by
$$d(K,L) = d(({\bf R}^n,||.||_K), ({\bf R}^n, ||.||_L)) .$$

Then it follows from Theorem 6 that
$$d(S(K,t),E(K^0)) \leq (\frac{6tv_n^2}{|K||K^0|})^{\frac{1}{2}}.$$
Thus for $\rho \in  {\bf R}$, $\rho > 1$,
$$d(\{x \in K: |K^x| \leq \rho |K^0| \},E(K^0)) \leq (6 \rho)^{\frac{1}{2}},$$
independent of $K$ and $n$.
It follows that for fixed $\rho$,  $\{x; |K^x|\leq \rho |K^0|\}$ is almost an
ellipsoid.
\newpage
\medskip\noindent{\bf {Proof of Theorem 6}}
\par
Let $u\in S^{n-1}, \lambda \in {\bf R}, 0 \leq \lambda <
\frac{1}{||u||_K}$ and $x=\lambda u$.
\newline
By (2) and symmetry
$$|K^x| = \int_0^{||u||_K} \phi(y)(\frac{1}{(1-\lambda y)^{n+1}}
+\frac{1}{(1+\lambda y)^{n+1}})dy,$$
where $\phi=\phi_{K^0}^u.$
For fixed $\lambda \geq0$ put
$$ f_\lambda(y) = \frac{1}{(1-\lambda y)^{n+1}}
+\frac{1}{(1+\lambda y)^{n+1}}.$$
Observe that $ f_\lambda$ is increasing in y, if $y \geq 0$. Put
$$a= \frac{n \int_0^\infty \phi(y)dy}{\phi(0)}=
\frac{n|K^0|}{2\phi(0)}.$$

Now we distinguish two cases.
\newline
1. $\lambda a < 1$.
\newline
Then we claim
that for all functions $\psi: {\bf R}_0^+
\rightarrow {\bf R}_0^+$ such that $\psi^{\frac{1}{n-1}}$ is continuous
on its support and continuous from the
right at 0, decreasing and concave on its support and such that
$$ (i) \hspace{.1in}\psi(0) = \phi(0)$$
$$ (ii) \int_0^\infty \psi(y)dy=\int_0^\infty \phi(y)dy=\frac{|K^0|}{2}$$
holds true, $\int_0^\infty \psi(y)f_\lambda(y)dy$ is maximal if $\psi$ is
of the
form
\[
\psi_0(y) = \left\{ \begin{array}{ll}
 \phi (0)(1-\frac{y}{a})^{n-1}  & \mbox{if $y\in [0,a]$}\\
0 & \mbox{otherwise.}
\end{array}
\right. \]

Indeed, let $\psi$ be a function with above properties and with support on
[0,$\tilde{a}$]. Put
$$H(t)= \int_{t}^{a} \psi(y)dy -\int_{t}^{a} \psi_0(y)dy.$$
Note that $\tilde{a} \leq a$, $H(0)=0 = H(a)$ and that the derivative of H
with respect to t is first negative, then positive; therefore $H(t) \leq 0$.
\newline
Consequently (with $g_\lambda (y) =f_\lambda (y)-2$)
$$\int_{0}^{\infty}
\psi(y)g_\lambda(y)dy=\int_{y=0}^{\infty}\psi(y)(\int_{t=0}^{y}
g^{'}_\lambda(t)dt)dy=
\int_{t=0}^{\infty}g^{'}_\lambda(t)(\int_{t}^{\infty}
\psi(y)dy)\hspace{.1in}dt$$
$$\leq\int_{0}^{\infty}g^{'}_\lambda(t)(\int_{t}^{\infty}
\psi_0(y)dy)\hspace{.1in}dt=\int_{0}^{\infty}
\psi_0(y)g_\lambda(y)dy.$$  From this the above claim follows.
\vskip 3mm
Hence
$$|K^x|\leq \phi(0)\int_0^a
(\frac{(1-\frac{y}{a})^{n-1}}{(1-\lambda y)^{n+1}}
+\frac{(1-\frac{y}{a})^{n-1}}{(1+\lambda y)^{n+1}})dy$$
$$=\frac{|K^0|}{1-\lambda^2 a^2}.$$
Here we have used Lemma 4 (iii).
\newline
For $x=\lambda u \in \partial S(K,t)$,
\begin{equation}
1=||x||_{S(K,t)}=\lambda ||u||_{S(K,t)}
\end{equation}
and
\begin{equation}
t=\frac{|K||K^x|}{v_n^2}.
\end{equation}
Therefore
$$t\leq \frac{|K||K^0|}{v_n^2(1-\lambda^2 a^2)}$$
and hence by (3)
$$||u||_{S(K,t)} \leq
\frac{n}{2}(1-\frac{|K||K^0|}{tv_n^2})^{-\frac{1}{2}}
\frac{|K^0|}{\phi(0)}.$$
Now (
[B] respectively [He]; see also [Mi-P])
\begin{equation}
\frac{|K^0|}{\phi(0)} \leq 2\sqrt3
(\frac{\int_{K^0}|<x,u>|^2dx}{|K^0|})^{\frac{1}{2}}=
2\sqrt3 ||u||_{E(K^0)},
\end{equation}
and thus
$$S(K,t)
\supseteq\frac{1}{\sqrt3\hspace{.1in}n}
(1-\frac{|K||K^0|}{tv_n^2})^{\frac{1}{2}}E(K^0).$$
\vskip 3mm
2. $\lambda a \geq 1$.
\newline
Again let $x=\lambda u \in \partial S(K,t)$. By definition of a, (3) and (5)
$$||u||_{S(K,t)} \leq\sqrt3\hspace{.1in} n||u||_{E(K^0)}.$$
This implies that
$$S(K,t)
\supseteq\frac{1}{\sqrt3\hspace{.1in}n}
E(K^0),$$
which proves the inclusion from below also in this case.
\vskip 3mm
On the other hand by
(1) and symmetry
$$|K^x|= \frac{1}{2n} \int_{S^{n-1}}
(\frac{1}{(||v||_{K^0} - <v,x>)^n}+
\frac{1}{(||v||_{K^0} + <v,x>)^n})d\sigma(v)$$
$$\geq \frac{1}{2n}\int_{S^{n-1}}
(2+n(n+1)(\frac{<x,v>}{||v||_{K^0}})^2)
\frac{d\sigma(v)}{||v||_{K^0}^n}$$
$$=|K^0| + \frac{(n+1)(n+2)}{2}\int_{K^0} |<x,y>|^2 dy$$
$$=|K^0| + \frac{(n+1)(n+2)}{2} \lambda^2|K^0| ||u||_{E(K^0)}^2.$$

Then (3) and (4) give
$$||u||_{S(K,t)} \geq
\frac{((n+1)(n+2))^{\frac{1}{2}}}{\sqrt2}
(\frac{tv_n^2}{|K||K^0|}-1)^{-\frac{1}{2}}
||u||_{E(K^0)}$$
or equivalently
$$S(K,t) \subseteq \frac{\sqrt2
}{((n+1)(n+2))^{\frac{1}{2}}}
(\frac{tv_n^2}{|K||K^0|}-1)^{\frac{1}{2}}E(K^0).$$
\vskip 5mm

Using (2) and a minimality argument similar to the
maximality argument of the above claim we get the other upper bound. Namely,
for fixed $\lambda$ and for
all functions $\psi: {\bf R}_0^+ \rightarrow {\bf R}_0^+$
such that $\psi^{\frac{1}{n-1}}$ is continuous on its support and
continuous from the right
at 0, decreasing and concave on its support and
such that
$$(i)\hspace{.1in} \psi(0) = \phi(0)$$
$$ (ii) \int_0^\infty \psi(y)dy=\int_0^\infty \phi(y)dy=\frac{|K^0|}{2}$$
holds true, $\int_0^\infty \psi(y)f_\lambda(y)dy$ is minimal if $\psi$ is
of the
form
\[
\psi(y) = \left\{ \begin{array}{ll}
 \phi(0)  & \mbox{if $y\in [0,a]$}\\
0 & \mbox{otherwise,}
\end{array}
\right. \]
where
$$a= \frac{|K^0|}{2\phi(0)}.$$
Note that in this situation $\lambda a < 1 $ always.
\newline
Consequently
$$|K^x|\geq \phi(0)\int_0^a f_\lambda(y) dy
\geq \frac{|K^0|}{(1-\lambda^2 a^2)^n}.$$
Then we use  again (3), (4) and the fact that ([B] respectively [He]; see
also [Mi-P])
$$\frac{|K^0|}{\phi(0)} \geq \sqrt2
||u||_{E(K^0)}$$
and get
$$S(K,t)\subseteq\sqrt2
(1-(\frac{|K||K^0|}{tv_n^2})^\frac{1}{n})^{\frac{1}{2}}E(K^0).$$

\vskip 5mm
Now we consider the non-symmetric case.
\vskip 3mm
\medskip\noindent{\bf Theorem 7}
\begin{it}
\par
Let K be a convex body in ${\bf R}^n$. Then
$$d_{n}^\prime(t) \hspace{.1in}E(K^0)
\subseteq S(K,t)
\subseteq
c_{n}^\prime(t) \hspace{.1in}E(K^0),$$
where
$$c_{n}^\prime(t)=\frac{2\sqrt2}{((e-2)(n+1)(n+2))^{\frac{1}{2}}}
(\frac{tv_n^2}{|K||K^0|})^{\frac{1}{2}}
(1-\frac{|K||K^0|}{tv_n^2})^{\frac{1}{2}}$$
and
$$d_{n}^\prime(t)=d_{n}(t)=\frac{1}{\sqrt3
\hspace{.1in}n}\hspace{.1in}(1-\frac{|K||K^0|}{tv_n^2})^{\frac{1}{2}}.$$
\end{it}
\vskip 5mm
\medskip\noindent{\bf {Proof}}
\par
By (2) we get for $u\in S^{n-1}$ and $x=\lambda u$ with $0 \leq
\lambda < \frac{1}{h_{K^0}(u)}$ that
$$|K^x|= \int_{-h_{K^0}(-u)}^{h_{K^0}(u)}\frac{\phi(y)}
{(1-\lambda y)^{n+1}}dy,$$
where $\phi=\phi_{K^0}^u$.
\newline
Notice that $K^0$ has its center of gravity at 0, as K
has its Santal${\mbox{\'o}}$-point at 0. Therefore
$$\int_{-h_{K^0}(-u)}^{h_{K^0}(u)}y\phi(y)dy=0.$$
Notice also that
\begin{equation}
1=h_{S(K,t)^0}(x)= \lambda h_{S(K,t)^0}(u).
\end{equation}
\vskip 3mm
Now we apply the following result of Fradelizi [F]
to the functions $\phi(y)$ and
$f_{\lambda}(y)=\frac{1}{(1-\lambda y)^{n+1}}$ to get the same
upper estimate for $|K^x|$as in the proof of Theorem 6.
Therefore $d_{n}^\prime(t)=d_{n}(t).$
\newpage
\medskip\noindent{\bf {Theorem}}([F])
\par
Let $\psi: {\bf R}
\rightarrow {\bf R}$, $\psi \geq 0$, such that $\psi^{\frac{1}{n-1}}$ is
continuous
and concave on its support and such that $\int_{-\infty}^\infty y\psi(y)dy =0$.
Let $f: {\bf R}
\rightarrow {\bf R}$ be any convex function. Then if
$a=\frac{n \int_{-\infty}^{\infty} \psi(y)dy}{2 \psi(0)}$, one has

$$ \int_{-\infty}^\infty \psi(y)f(y)dy \leq
\psi(0) \int_{-a}^a (1-\frac{|y|}{a})^{n-1}f(y)dy.$$

\vskip 7mm
For the right-hand side inclusion write
$$|K^x|=\int_{-h_{K^0}(-u)}^{0}\frac{\phi(y)}
{(1-\lambda y)^{n+1}}dy + \int_{0}^{h_{K^0}(u)}\frac{\phi(y)}
{(1-\lambda y)^{n+1}}dy$$
$$\geq$$
$$\int_{-h_{K^0}(-u)}^{0}\phi(y)(1+(n+1)\lambda y)dy
+\int_{0}^{h_{K^0}(u)}\phi(y)(1+(n+1)\lambda y +
\frac{(n+1)(n+2)}{2}\lambda^2 y^2)dy$$
$$=|K^0| + \frac{(n+1)(n+2)}{2}\lambda^2\int_{0}^{h_{K^0}(u)}y^2\phi(y)dy,$$

where for the last equality we have used the fact that the center of
gravity is at $0$.
\newline
Let a be such that
$$\int_{0}^{h_{K^0}(u)}y\phi(y)dy=\phi(0)\frac{a^2}{n(n+1)}=
\int_{-a}^{\frac{a}{n}}y\psi_0(y)dy,$$
where $\psi_0(y) =\phi(0)(1+\frac{y}{a})^{n-1}.$
\par
Now one shows as in the beginning of the proof of Theorem 6 that
for all functions $\psi: {\bf R}_0^+
\rightarrow {\bf R}_0^+$ such that $\psi^{\frac{1}{n-1}}$ is continuous
on its support and continuous from the
right at 0, concave on its support and such that
$$ (i) \hspace{.1in}\psi(0) = \phi(0)$$
$$ (ii) \int_0^\infty y\psi(y)dy=\int_0^\infty y\phi(y)dy$$
holds true, $\int_0^\infty y^2\psi(y)dy$ is minimal if $\psi$ is of the
form
\[
\psi_0(y) = \left\{ \begin{array}{ll}
 \phi (0)(1+\frac{y}{a})^{n-1}  & \mbox{if $y\in [0,\frac{a}{n}]$}\\
0 & \mbox{otherwise}
\end{array}
\right  .\]
and that for all functions $\psi: {\bf R}_0^-
\rightarrow {\bf R}_0^+$ such that $\psi^{\frac{1}{n-1}}$ is continuous
on its support and continuous from the
left at 0, concave on its support and such that
$$ (i) \hspace{.1in}\psi(0) = \phi(0)$$
$$ (ii) \int_{-\infty}^0 y\psi(y)dy=\int_{-\infty}^0 y\phi(y)dy$$
holds true, $\int_{-\infty}^0 y^2\psi(y)dy$ is maximal if $\psi$ is of the
form
\[
\psi_0(y) = \left\{ \begin{array}{ll}
 \phi (0)(1+\frac{y}{a})^{n-1}  & \mbox{if $y\in [-a,0]$}\\
0 & \mbox{otherwise.}
\end{array}
\right. \]

Therefore
$$\int_{0}^{h_{K^0}(u)}y^2\phi(y)dy \geq \int_{0}^{\frac{a}{n}}y^2\psi_0(y)dy =
\frac{(1+\frac{1}{n})^{n+1} -2}{n(n+1)(n+2)}a^3 \phi(0)$$
and
$$\int_{-h_{K^0}(-u)}^0 y^2\phi(y)dy \leq \int_{-a}^{0}y^2\psi_0(y)dy =
\frac{2}{n(n+1)(n+2)}a^3 \phi(0).$$
Thus we get
$$\int_{0}^{h_{K^0}(u)}y^2\phi(y)dy \geq
\frac{e-2}{2}\int_{-h_{K^0}(-u)}^0 y^2\phi(y)dy,$$
and therefore
$$\int_{0}^{h_{K^0}(u)}y^2\phi(y)dy \geq
\frac{e-2}{4}\int_{-h_{K^0}(-u)}^{h_{K^0}(u)} y^2\phi(y)dy.$$
It follows that
$$|K^x| \geq |K^0| (1+ \frac{(n+1)(n+2)(e-2)}{8}
\lambda ^2 \|u\|_{E(K^0)}^2),$$
which implies, using (6),
$$S(K,t) \subseteq \frac{2\sqrt2}{((e-2)(n+1)(n+2))^{\frac{1}{2}}}
(\frac{tv_n^2}{|K||K^0|})^{\frac{1}{2}}
(1-\frac{|K||K^0|}{tv_n^2})^{\frac{1}{2}} E(K^0).$$
\vskip 7mm
Next we give estimates on the ``size" of $S(K,t)$ in terms of the body
K.
\newline
We also need the following Lemma.

\vskip 3mm
\medskip\noindent{\bf Lemma 8} (see for instance [S2])
\par
\begin{it}
Let K be a convex body in ${\bf R}^n$ such that the center of gravity of
K is at 0. Then

$$\frac{1}{e} \leq \frac{1}{|K|}\hspace{.1in} \int_0^{h_{K}(u)} \phi(y)dy
\hspace{.1in}\leq 1-\frac{1}{e}.$$
\end{it}
\newpage

\medskip\noindent{\bf Theorem 9}
\par
\begin{it}
Let K be a convex body in ${\bf R}^n$. Then
\newline
(i)
$$(1-(\frac{|K||K^0|}{tv_n^2})^{\frac{1}{n}})\hspace{.1in}K
\subseteq \hspace{.1in}S(K,t)\hspace{.1in}
\subseteq (1-\frac{|K||K^0|}{etv_n^2})\hspace{.1in}K.$$
(ii) If in addition K is symmetric, then
$$(1-(\frac{|K||K^0|}{tv_n^2})^{\frac{1}{n}})\hspace{.1in}K
\subseteq \hspace{.1in}S(K,t)\hspace{.1in}
\subseteq (1-\frac{|K||K^0|}{tv_n^2})^{\frac{1}{2}}\hspace{.1in}K.$$
\end{it}
\vskip 3mm
\medskip\noindent{\bf {Proof}}
\par
Let $u \in S^{n-1}$, $\lambda \in {\bf R}, 0 \leq \lambda \leq1,$ be given
and let
$x=\frac{\lambda}{h_{K^0}(u)}u$. Then $K$ contains $\alpha K + x$
for all $\alpha$, $0 \leq \alpha \leq 1-\lambda$ and consequently $K^x
\subseteq
\frac{1}{\alpha} K^0$; therefore we have for $x \in \partial S(K,t)$
$$ t=\frac{|K||K^x|}{v_n^2}\leq \frac{|K||K^0|}{\alpha^n v_n^2}$$
and hence
$$\alpha \leq (\frac{|K||K_0|}{tv_n^2})^{\frac{1}{n}}.$$
Thus we have for all $\lambda$ with $\lambda
\leq 1-(\frac{|K||K_0|}{tv_n^2})^{\frac{1}{n}}$ that
$$\frac{\lambda}{h_{K^0}(u)}u \in S(K,t).$$
This proves the left-hand side inclusion.
\vskip 3mm
For the right-hand side inclusion we first treat the
symmetric case.
\newline
Let $x=\lambda u, u \in S^{n-1}$, $0 \leq  \lambda < ||u||_K^{-1}$.
Let $f_\lambda$ be as in the proof of Theorem 6. By (2) and symmetry
$$|K^x| = \int_0^{||u||_K} \phi(y)f_\lambda(y) dy.$$
Notice that for all functions $\psi: {\bf R}_0^+ \rightarrow {\bf R}_0^+$
such that $\psi^{\frac{1}{n-1}}$ is continuous on its support and continuous
at 0 from the right, decreasing and  concave on its support and such that
$$ (i)\hspace{.1in}\psi > 0 \hspace{.1in}\mbox{on}\hspace{.1in}
[0,||u||_K),
\hspace{.2in}\psi =
0\hspace{.1in}\mbox{on} \hspace{.1in}[||u||_K,
\infty)$$
$$ (ii) \int_0^{||u||_K} \psi(y)dy=\int_0^{||u||_K}
\phi(y)dy=\frac{|K^0|}{2}$$ holds true, $\int_0^{||u||_K}
\psi(y)f_\lambda(y)dy$
is smallest if $\psi$ is of the
form
\[
\psi(y) = \left\{ \begin{array}{ll}
 c(1-\frac{y}{{||u||_K}})^{n-1}  & \mbox{if $y\in [0,{||u||_K})$}\\
0 & \mbox{otherwise,}
\end{array}
\right. \]
where
$$c= \frac{n|K^0|}{2{||u||_K}}.$$
Hence
$$|K^x|\geq c\int_0^{||u||_K}
(\frac{(1-\frac{y}{||u||_K})^{n-1}}{(1-\lambda y)^{n+1}}
+\frac{(1-\frac{y}{||u||_K})^{n-1}}{(1+\lambda y)^{n+1}})dy$$
$$=\frac{|K^0|}{1-\lambda^2 ||u||_K^2},$$
which implies
$$S(K,t) \subseteq (1-\frac{|K||K_0|}{tv_n^2})^{\frac{1}{2}}K.$$
\vskip 3mm
Next we consider the non-symmetric case.
$$|K^x|= \int_0^{h_{K^0}(u)}\frac{\phi(y)}{(1-\lambda y)^{n+1}}dy +
\int_0^{h_{K^0}(-u)}\frac{\phi(-y)}{(1+\lambda y)^{n+1}}dy$$
$$\geq\int_0^{h_{K^0}(u)}\frac{\phi(y)}{(1-\lambda y)^{n+1}}dy .$$

Fix $\lambda$ and note again that among all functions
$\psi:{\bf R}_0^+ \rightarrow {\bf R}_0^+$\hspace{.1in}such that
$\psi^\frac{1}{n-1}$ is continuous on its support and continuous
from the right
at 0, decreasing and concave  on its support and for which

$$ \psi > 0 \hspace{.1in}\mbox{on}\hspace{.1in} [0,h_{K^0}(u)),
\hspace{.2in}\psi = 0\hspace{.1in}\mbox{on} \hspace{.1in}[h_{K^0}(u),
\infty),$$

$$ \int_0^{h_{K^0}(u)} \psi(y)dy=\int_0^{h_{K^0}(u)}
\phi(y)dy$$

\vskip 3mm
holds true, $\int_0^{h_{K^0}(u)}\frac{\psi(y)}{(1-\lambda y)^{n+1}}dy$
is smallest if  $\psi$ is of the
form
\[
\psi(y) = \left\{ \begin{array}{ll}
 c(1-\frac{y}{h_{K^0}(u)})^{n-1}  & \mbox{if $y\in [0,h_{K^0}(u))$}\\
0 & \mbox{otherwise,}
\end{array}
\right. \]
where
$$c= \frac{n \int_0^{h_{K^0}(u)} \phi(y)dy}{h_{K^0}(u)}.$$

Arguments similar to the ones before together with Lemma 8 then
finish the proof.

\newpage
II. SANTAL${\mbox{\'O}}$-REGIONS AND AFFINE SURFACE AREA
\vskip 3mm
Recall that for a convex body K in ${\bf R}^n$ the affine surface area
$$as (K) = \int_{\partial K} \kappa(x)^{\frac{1}{n+1}} d\mu(x),$$
where $\kappa(x)$ is the (generalized) Gaussian curvature in
$x \in \partial K$ and $\mu$ is the surface measure on $\partial K$.
We prove here
\vskip 5mm
\medskip\noindent {\bf Theorem 10}
\par
\begin{it}
Let K be a convex body in ${\bf R}^n$. Then
$$\mbox{lim}_{t\rightarrow \infty} t^{\frac{2}{n+1}}(|K| -|S(K,t)|)
=
\frac{1}{2}
(\frac{|K|}{v_n })^{\frac{2}{n+1}}\hspace{.1in} as(K).$$

\end{it}

\vskip 5mm
In the proof of Theorem 10 we follow the ideas of [S-W].
We need several Lemmas for the proof. We also use the following notations.
For $x \in
\partial K$,
$N(x)$ is
the  outer unit normal vector to $\partial K$ in $x$. For two points x and y in
${\bf R}^n$, $[x,y]=\{\alpha x +(1-\alpha )y: 0 \leq \alpha \leq 1\}$
denotes the
line segment from x to y.
\vskip 5mm
The proof of the following Lemma is standard.
\newline
\medskip\noindent {\bf Lemma 11}
\par

\begin{it}
Let K and L be two convex bodies in ${\bf R}^n$ such that $0 \in
\mbox{int} (L)$ and $L \subseteq K$.
Then
$$|K| - |L| = \frac{1}{n} \int_{\partial K}
<x, N(x)>(1-(\frac{||x_L||}{||x||})^n) d\mu(x),$$
where $x_L = [0,x] \cap \partial L$.
\end{it}

\vskip 5mm
For $x \in \partial K$ denote by $r(x)$ the radius of the  biggest
Euclidean ball contained in K that touches $\partial K$ at x.
More precisely
$$r(x)=\mbox{max}\{r:x\in B(y,r) \subset K \hspace{.1in}\mbox{for some y}
\in K\}.$$
\vskip 5mm

\medskip\noindent {\bf Remark}

It was shown in [S-W] that

(i) If $B(0,1) \subset K$, then
$$\mu \{ x \in \partial K: r(x) \geq \beta \} \geq (1-\beta)^{n-1}
\mbox{vol}_{n-1} (\partial K)$$

 (ii) $$\int_{\partial K} r(x)^{- \alpha} d \mu (x) < \infty \qquad
\mbox{for all} \quad \alpha ,
\quad 0 < \alpha < 1$$
\vskip 5mm
We postpone the proof of the next two lemmas which we use for the proof of
Theorem 10.
\vskip 3mm

\medskip\noindent {\bf Lemma 12}
\begin{it}
Suppose $0$ is in the interior of $K$.  Then we have for all $x$ with
$r(x) > 0$ and for all t such that $(S(K,t)$ has non-empty interior.
$$0 \leq \frac{1}{n} < x,N(x) > t^{\frac{2}{n+1}} \left(1- (
\frac{\| x_t \|}{\| x \|} )^n \right) \leq c\hspace{.1in}r(x)^{-
\frac{n-1}{n+1}},$$ where $x_t=[0,x]\cap \partial S(K,t)$ and
$c$ is a constant
independent  of $x$ and $t$.

\end{it}

\vskip 5mm
\medskip\noindent {\bf Lemma 13}
\begin{it}
Suppose $0$ is in the interior of $K$. Then
$$\lim_{t \rightarrow \infty} \frac{1}{n} < x,N(x) >
t^{\frac{2}{n+1}} \left(1- ( \frac{\| x_t \|}{\| x \|} )^n
\right)$$
exists a.e.
and is equal to $$\frac{1}{2}(\frac{|K|}{v_n})^{\frac{2}{n+1}}
\kappa(x)^{\frac{1}{n+1}},$$
where $\kappa (x)$ is the Gaussian curvature in $x \in \partial K$.
\end{it}

\vskip 10mm
\medskip\noindent {\bf Proof of Theorem 10}
\vskip 3mm
We may assume that $0$ is in the interior of $K$.  By Lemma 11 and with the
notations of Lemma 12 we have

$$|K| - |S(K,t)| = \frac{1}{n} \int_{\partial K}
<x, N(x)>(1-(\frac{||x_t||}{||x||})^n) d\mu(x)$$
By Lemma 12 and the Remark preceding it, the functions under the
integral sign are
bounded uniformly in $t$ by an $L^1$-function and by Lemma 13 they are
converging pointwise a.e.
We apply Lebesgue's convergence theorem.
\newpage

\medskip\noindent {\bf {Proof of Lemma 12}}
\par
Let $x \in \partial K$ such that $r(x) > 0$.
As
$\| x_t\|=\| x\|- \| x-x_t \|$, we have

\begin{equation}
\frac{1}{n} <x,N(x) > \left(1- ( \frac{ \|
x_t
\|}{\| x \|} )^n  \right)
\leq \hspace{.1in}<\frac{x}{\| x \|}, N(x)> \| x \ - x_t\|.
\end{equation}
\vskip 3mm
a) We consider first the case where
$$\|x-x_t\|\hspace{.1in} < \hspace{.1in}r(x) <\frac{x}{\|x\|},
N(x)>.$$
Let $\tilde{\rho} =
\|x_t -(x-r(x) N(x))\|$. By
assumption
$0 <
\tilde{\rho} < r(x)$. Computing $\tilde{\rho}$ we get
$$\tilde{\rho}= (||x-x_t||^2 +r(x)^2
-2r(x)||x-x_t||<\frac{x}{\|x\|},N(x)>)^{1/2}.$$
Since K contains the Euclidean
ball of radius r(x) centered at
$x-r(x)N(x)$, $K^{x_t}$ is contained in the polar (with
respect to $x_t$) of the Euclidean ball with radius r(x). Hence by
Corollary 5,
$$ t= \frac{|K||K^{x_t}|}{v_n^2} \leq
\frac{|K|}{v_nr(x)^n (1-(\frac{\tilde{\rho}}{r(x)})^2)^{\frac{n+1}{2}}}$$
and therefore, using (7),
$$\frac{1}{n} < x,N(x) >
t^{\frac{2}{n+1}} (1- ( \frac{\| x_t \|}{\| x \|} )^n )
\leq
(\frac{|K|}{v_n})^{\frac{2}{n+1}} r(x)^{-\frac{n-1}{n+1}},
$$
which proves Lemma 12 in this case.
\vskip 3mm
b) Now we consider the case where $$\|x-x_t\|
\geq \hspace{.1in} r(x)\hspace{.1in} <\frac{x}{\|x\|},
N(x)>.$$
We can suppose that $t$ is big enough so that $x_t \neq 0$.
We choose $\alpha > 0$ such that
$
B(0,\alpha) \subseteq K \subseteq
B(0,\frac{1}{\alpha})$ and t so big that $x_t \notin B(0,\alpha)$.
K contains the spherical cone
$C=\mbox{co}[x,H\cap B(0,\alpha)]$,
where H is the hyperplane through 0 orthogonal to the line segment [0,x].
We get
$$|C^{x_t}| = \frac{v_{n-1} ||x||^n}{n \alpha^{n-1}
||x_t||(||x||-||x_t||)^n}.$$
Consequently
$$t = \frac{|K||K^{x_t}|}{v_n^2}
\leq
\frac{|K|v_{n-1}\| x\|^n}
{n  v_n^2\alpha^{n-1}
\| x_t\|\quad\|x-x_t\|^n}$$
and hence, using (7),
$$\frac{1}{n} < x,N(x) >
t^{\frac{2}{n+1}} (1- ( \frac{\| x_t \|}{\| x \|} )^n )
\leq
(\frac{|K| v_{n-1}}{n\quad v_n^2})^{\frac{2}{n+1}}
r(x)^{-\frac{n-1}{n+1}}\frac{1}{\alpha^{\frac{4n}{n+1}}}.$$

\vskip 5mm
\medskip\noindent {\bf {Proof of Lemma 13}}
\par
As in the proof of Lemma 12, we can choose an $\alpha > 0$ such
that

$$B(0,\alpha) \subseteq K \subseteq
B(0,\frac{1}{\alpha}).$$ Therefore
\begin{equation}
1 \geq \hspace{.1in}< \frac{x}{\| x \|} ,N(x) >\hspace{.1in} \geq \alpha^2.
\end{equation}

Since $x$ and
$x_t$ are colinear,
$$||x||=||x_t||+||x-x_t||,$$
and hence
\begin{equation}
\frac{1}{n} < x,N(x) > \left(1- ( \frac{\| x_t\|}{\|
x
\|} )^n \right)=\frac{1}{n} < x,N(x) > \left(
(1-(1-
\frac{\| x -x_t \|}{\| x \|} )^n \right)$$
$$\geq$$
$$< \frac{x}{\| x \|} ,N(x) > \| x -x_t \| \left( 1 -
d \cdot
\frac{\| x -x_t \|}{\| x \|} \right),
\end{equation}
for some constant $d$, if we choose $t$ sufficiently large.
We denote by $\theta$ the angle between $x$ and $N(x)$. Then
$< \frac{x}{\| x \|} ,N(x) > = \mbox{cos} \theta $.
\vskip 3mm
By [S-W] $r(x) > 0$ a.e.
and by [L2] the indicatrix of Dupin exists a.e. and is an
elliptic cylinder or an ellipsoid.
\vskip 3mm
(i) \underline{Case where the indicatrix is an ellipsoid}
\newline
This case can be reduced
to the case
of a sphere by
an affine transformation with determinant 1 (see for instance  [S-W]). Let
$\sqrt{\rho(x)}$ be the radius of this sphere.
Recall that we have to show that
$$\lim_{t \rightarrow \infty} \frac{1}{n} < x,N(x) >
t^{\frac{2}{n+1}} \left(1- ( \frac{\| x_t \|}{\| x \|} )^n
\right)
=\frac{1}{2}(\frac{|K|}{v_n})^{\frac{2}{n+1}} \rho
(x)^{-\frac{n-1}{n+1}}.$$
We put $\rho(x)=\rho$ and we introduce a coordinate system such that $x =0$ and
$N(x) = (0, \ldots 0, -1)$.  $H_0$ is the tangent hyperplane to $\partial K$ in
$x =0$ and $\{ H_s: s \geq 0 \}$ is the family of hyperplanes
parallel to $H_0$
that have non-empty
intersection with $K$ and are at distance $s$ from $H_0$.  For $s > 0$,
$H_s^+$ is the
halfspace generated by $H_s$ that contains $x =0$.
For $a\in {\bf R}$, let $z_a=(0,\ldots 0,a)$ and $B_a=B(z_a,a)$ be the
Euclidean ball
with center $z_a$ and radius a.
As in [W], for
$\epsilon > 0$ we can choose $s_0$ so small that for all $s \leq s_0$

$$B_{\rho-\varepsilon} \cap
H_s^+
\subseteq
K \cap H_s^+
\subseteq
B_{\rho+\varepsilon}
\cap H_s^+.$$
For $\lambda \in {\bf R}$ let
$G_{\lambda}=\{ x: <x,z_{\rho+\varepsilon}-x_t>=\lambda\}$ be a
hyperplane orthogonal to the line segment $[x_t,z_{\rho+\varepsilon}]$, if
$t$ is sufficiently large.
Let $\lambda_0=\mbox{max} \{\lambda: G_{\lambda}^+\cap B_{\rho+\varepsilon}
\subseteq  H_{s_0}^+ \cap B_{\rho+\varepsilon} \}$.
Define C to be the cone tangent to $B_{\rho+\varepsilon}$ at
$G_{\lambda_0} \cap B_{\rho+\varepsilon}$ and choose the minimal $\lambda_1$
so that
$$K \cap \{x: \lambda_0 \leq <x,z_{\rho+\varepsilon}-x_t> \leq \lambda_1\}
\subseteq D=C \cap \{x: \lambda_0 \leq <x,z_{\rho+\varepsilon}-x_t> \leq
\lambda_1\}.$$
Then K is contained in the union of the truncated cone D of
height $h=|\lambda_1-\lambda_0|$ and the cap
$L=\{x \in B_{\rho +\varepsilon}:<x,z_{\rho+\varepsilon}-x_t> \leq
\lambda_0 \}$
(see Figure 1).
\begin{figure}
\BoxedEPSF{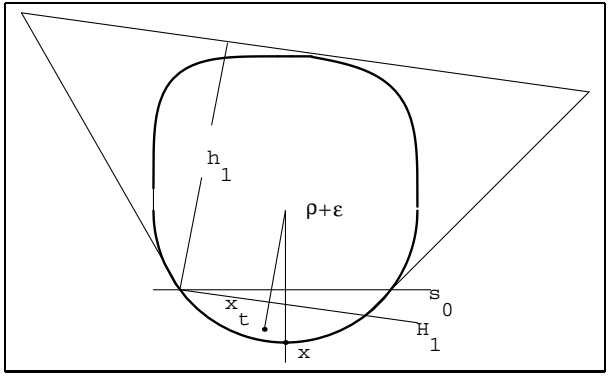 scaled 800}
\caption{the estimate from below}
\end{figure}

Therefore
$$ K^{x_t}\supseteq (D \cup L)^{x_t},$$

and to estimate $ |K^{x_t}|$ we have to compute $|(D \cup L)^{x_t}|$.
To do so we prove the following more general result.
\vskip 3mm
\underline{Claim 1}
\newline

\begin{it}
Let M be the convex body that is  the union of a truncated spherical
cone D with height h and a cap L  of a Euclidean ball with radius r such that
D is ``tangent" to L. For a point x in L and on the axis of symmetry of M
let  a=distance(x,D), b=distance(x,$\partial$L) and
$b_0=\frac{rb +(a+b)(r-b)}{r-(a+b)}$
(see Figure 2).
Then, if $x$ is such that $r > a+b$,
$$|M^x|=v_{n-1}(\frac{1}{r^n}\int_{\frac{r}{(r-b)+b_0}}
^1\frac{(1-y^2)^{\frac{n-1}{2}}}{(1-\frac{(r-b) y}{r})^{n+1}}dy$$
$$+\frac{1}{n}(\frac{1}{b_0}+
\frac{1}{a+h})(\frac{(2r(a+b)-(a+b)^2)^{\frac{1}{2}}}{rb
+(a+b)(r-b)})^{n-1}).$$
\end{it}
\vskip 3mm
\underline{Proof of Claim 1}
\newline
We introduce a coordinate system such that x=0 and such that
the $x_1$-axis coincides with the axis of symmetry of M (see Figure 2).
\begin{figure}
\BoxedEPSF{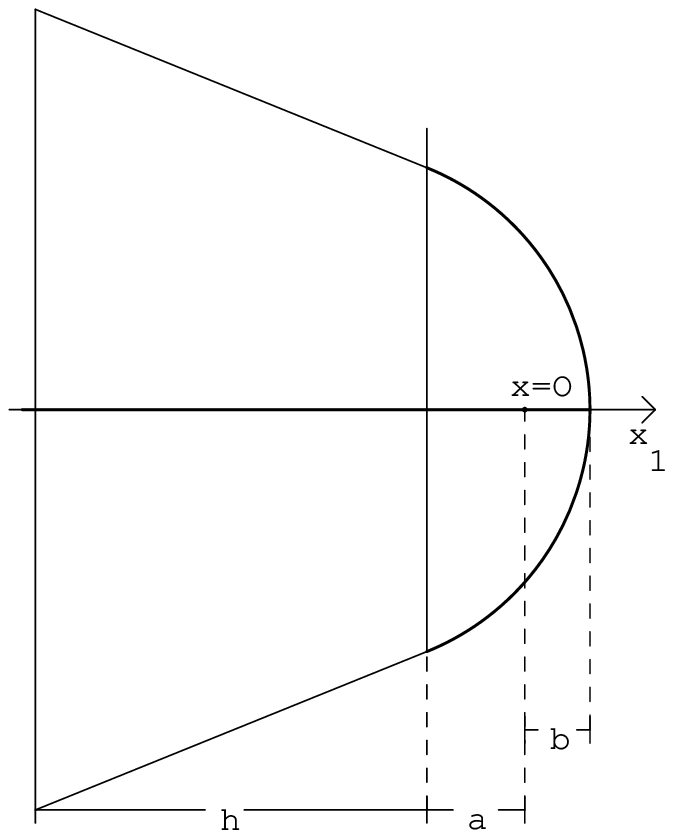 scaled 800}
\caption{Claim 1}
\end{figure}

Notice now that $M^0$ is such that each (n-1)-dimensional section
orthogonal to the $x_1$-axis is an (n-1)-dimensional Euclidean ball with
radius
$l(x_1)$, where
$$ l(x_1) = \frac{((a+h)x_1+1)(2r(a+b)-(a+b)^2)^{\frac{1}{2}}}
{(a+h+b_0)(r-(a+b))},\hspace{.2in}
\mbox{if}\hspace{.1in}-\frac{1}{a+h} \leq x_1
\leq \frac{1}{b_0},$$
$$ l(x_1) = \frac{1}{r}((1+x_1(r-b))^2-r^2x_1^2)^\frac{1}{2},
\hspace{.2in}\mbox{if}\hspace{.2in}\frac{1}{b_0} \leq x_1
\leq \frac{1}{b}.$$
From this Claim 1 follows.
\vskip 5mm
Now we apply Claim 1 to our situation. Then
$$r=\rho+\epsilon\hspace{.1in} \mbox{and} \hspace{.1in}
b=\rho+\epsilon-c,$$ where
$$c^2= ||x-x_t||^2 +(\rho+\epsilon)^2
- 2(\rho+\epsilon)||x-x_t||\mbox{cos}\theta.$$
$$a=c-\frac{\rho+\epsilon-s_0}{c}(c^2-||x-x_t||^2
\mbox{sin}^2 \theta)^{\frac{1}{2}}-
\frac{(2(\rho+\epsilon)s_0-s_0^2)^{\frac{1}{2}}||x-x_t||
\mbox{sin}\theta}{c}$$

and
$$b_0=\frac{(\rho+\epsilon)^2}{c-a}-c.$$
Therefore
$$|K^{x_t}| \geq
\frac{v_{n-1}}{(\rho+\epsilon)^n}\int_0
^1\frac{(1-y^2)^{\frac{n-1}{2}}}{(1-\frac{c
y}{\rho+\epsilon})^{n+1}}dy$$
$$+\frac{1}{n}(\frac{1}{b_0}+
\frac{1}{a+h})\frac{v_{n-1}}{r_{D}^{n-1}}-
\frac{v_{n-1}}{(\rho+\epsilon)^n}\int_0
^{\frac{\rho+\epsilon}{b_0+c}}\frac{(1-y^2)^{\frac{n-1}{2}}}
{(1-\frac{c
y}{\rho+\epsilon})^{n+1}}dy,$$
where $r_{D}$ is the radius of the base of the spherical
cone in $(D\cup L) ^{x_t}$. We put
$$R=\frac{1}{n}(\frac{1}{b_0}+
\frac{1}{a+h})\frac{v_{n-1}}{r_{D}^{n-1}}-
\frac{v_{n-1}}{(\rho+\epsilon)^n}\int_0
^{\frac{\rho+\epsilon}{b_0+c}}\frac{(1-y^2)^{\frac{n-1}{2}}}
{(1-\frac{c
y}{\rho+\epsilon})^{n+1}}dy.$$
Then, by Lemma 4 (ii), for $\varepsilon > 0$
$$|K^{x_t}| \geq
\frac{(1-\epsilon)v_n}{2^{\frac{n+1}{2}}(\rho+\epsilon)^{n}
(\frac{c}{\rho+\epsilon})^{\frac{n+1}{2}}
(1-\frac{c}{\rho+\epsilon})^{\frac{n+1}{2}}} +R,$$
provided that
$$\frac{c}{\rho+\epsilon} > \frac{1}{1+\frac{2\epsilon}{3n}}.$$
We choose t so big that this holds.
\newline
Hence
$$t=\frac{|K||K^{x_t}|}{v_n^2} \geq$$
$$(\frac{1}{2})^{\frac{n+1}{2}}(\frac{(1-\epsilon)|K|}{v_n})
\frac{(\rho+\epsilon)^{-n}}{(\frac{c}{\rho+\epsilon}
(1-\frac{c}{\rho+\epsilon}))^{\frac{n+1}{2}}}
\{1+
\frac{2^{\frac{n+1}{2}}R(\rho+\epsilon)^n(\frac{c}
{\rho+\epsilon})^{\frac{n+1}{2}}
(1-\frac{c}{\rho+\epsilon})^{\frac{n+1}{2}}}
{(1-\epsilon)v_n}\}$$
and by (9) for some constant d
$$\frac{t^{\frac{2}{n+1}}}{n} < x,N(x) > \left(1- ( \frac{\|
x_t\|}{\| x
\|} )^n \right)$$
$$\geq$$
$$\frac{1}{2}
(\frac{(1-\epsilon)|K|}{v_n})^{\frac{2}{n+1}}
(\rho+\epsilon)^{-\frac{n-1}{n+1}}
\frac{(1+k(2\frac{\| x -x_t\|\mbox{cos}\theta}
{\rho+\epsilon}-\frac{\| x -x_t
\|^2}{(\rho+\epsilon)^2}))^{-1}
(1 -d \cdot\frac{\| x -x_t \|}{\| x \|})}{(1+\frac{\| x -x_t
\|^2}{(\rho+\epsilon)^2}-2\frac{\| x -x_t\|\mbox{cos}\theta}
{\rho+\epsilon})^{\frac{1}{2}}(1-\frac{\| x -x_t\|}
{2(\rho+\epsilon)\mbox{cos}\theta})}$$
$$\cdot\{1+\frac{2^{\frac{n+1}{2}}R(\rho+\epsilon)^n}
{(1-\epsilon)v_n}(1-\frac{c}{\rho+\epsilon})^{\frac{n+1}{2}}
(\frac{c}{\rho+\epsilon})^{\frac{n+1}{2}}
\}^{\frac{2}{n+1}},$$
as
$$1-\frac{c}{\rho+\epsilon}\leq \frac{\| x -x_t\|\mbox{cos}\theta}
{\rho+\epsilon}(1-\frac{\| x -x_t\|}
{2(\rho+\epsilon)\mbox{cos}\theta})(1+k(2\frac{\| x  x_t\|\mbox{cos}\theta}
{\rho+\epsilon}-\frac{\| x -x_t
\|^2}{(\rho+\epsilon)^2})),$$
for some constant k, if t is big enough. $R$ remains bounded for $t
\rightarrow \infty$.
\newline
Note also that $\mbox{cos}\theta$ remains bounded from below  by (8).
\newline

Thus we have a lower bound for the expression in question.
\vskip 7mm
To get an upper bound we proceed in a similar way.
For $\lambda \in {\bf R}$ let now
$G_{\lambda}=\{ x: <x,z_{\rho-\varepsilon}-x_t>=\lambda\}$ be a
hyperplane orthogonal to the line segment $[x_t,z_{\rho-\varepsilon}]$.
Let $\lambda_0=\mbox{max} \{\lambda: G_{\lambda}^+\cap B_{\rho \varepsilon}
\subseteq  H_{s_0}^+ \cap B_{\rho-\varepsilon} \}$.
Let P be the point where the half-line starting at $x_t$ through
$z_{\rho -\epsilon}$ intersects $\partial K$.

Let C be the spherical cone $C=\mbox{co}[P, B_{\rho-\varepsilon} \cap
G_{\lambda_0}].$
Let h be the height of this cone (see Figure 3).

\begin{figure}
\BoxedEPSF{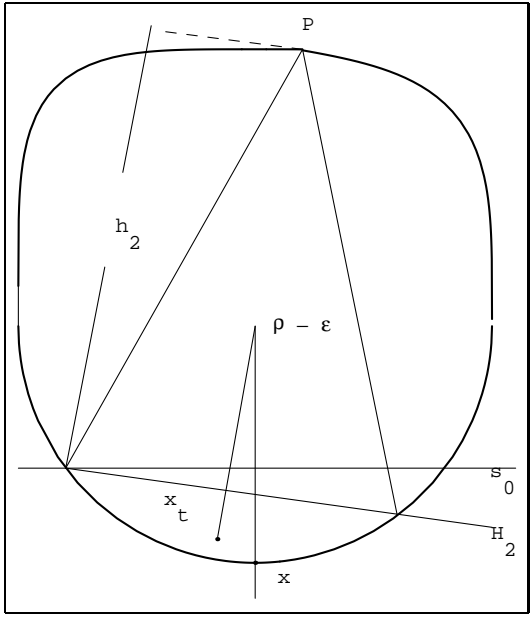 scaled 800}
\caption{the estimate from above}
\end{figure}

Let
$L=\{x \in B_{\rho -\varepsilon}:<x,z_{\rho-\varepsilon}-x_t> \leq
\lambda_0 \}$.
Then $K \supseteq C \cup L$,  and hence
$$K^{x_t} \subseteq (C \cup L)^{x_t},$$

and to estimate $ |K^{x_t}|$ we have to compute
$(C \cup L)^{x_t}$.
\vskip 3mm
To do so we prove the more general
\vskip 3mm
\underline{Claim 2}
\newline
\begin{it}
Let M be the union of a spherical cone C with height h and a
cap L  of a Euclidean ball with radius r such that the base of C
coincides with  the base of L. For a point x in L and on the axis of
symmetry of M
let  $\alpha=distance(x,C)$, $\beta=distance(x,\partial L)$ (see Figure 4).
Then with $\beta_0=\frac{r\beta +(\alpha+\beta)(r-\beta)}{r
(\alpha+\beta)}$ and
$x$ chosen such that $r> \alpha + \beta$,
$$|M^x|=v_{n-1}(\frac{1}{r^n}\int_{\frac{r}{(r-\beta)+\beta_0}}
^1\frac{(1-y^2)^{\frac{n-1}{2}}}{(1-\frac{(r-\beta) y}{r})^{n+1}}dy$$
$$+\frac{1}{n \alpha
(2r(\alpha+\beta)-(\alpha+\beta)^2)^{\frac{n
1}{2}}}(\frac{(\alpha+\beta_0)^n}{\beta_0^n}-
\frac{h^n}{(\alpha+h)^n})).$$
\end{it}
\newpage
\underline{Proof of Claim 2}
\newline
We introduce a coordinate system such that x=0 and such that
the $x_1$-axis coincides with the axis of symmetry of M (see Figure 4).
\begin{figure}
\BoxedEPSF{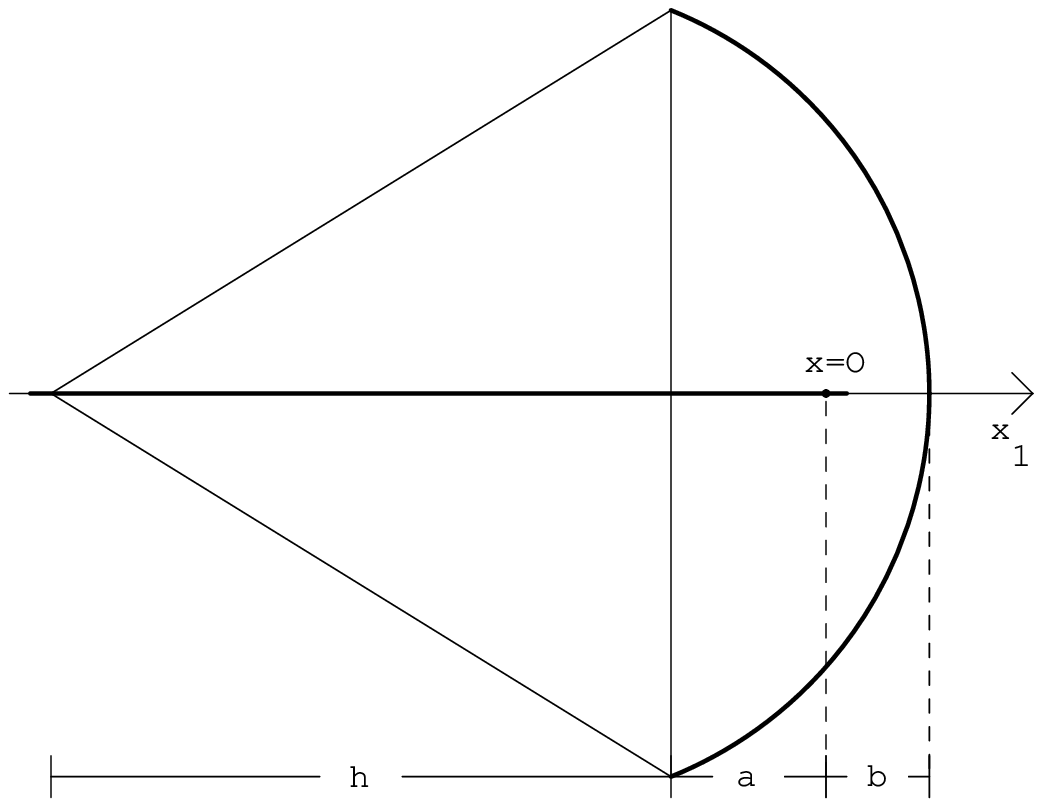 scaled 800}
\caption{Claim 2}
\end{figure}

Notice now that $M^0$ is such that each (n-1)-dimensional section
orthogonal to the $x_1$-axis is an (n-1)-dimensional Euclidean ball with
radius
$l(x_1)$, where
$$ l(x_1) = \frac{\alpha x_1+1}{(2r(\alpha +\beta)-(\alpha
+\beta)^2)^{\frac{1}{2}}},
\hspace{.2in}\mbox{if}\hspace{.1in}-\frac{1}{\alpha +h} \leq x_1
\leq \frac{1}{\beta_0},$$
$$ l(x_1) = \frac{1}{r}((1+x_1(r-\beta))^2-r^2x_1^2)^\frac{1}{2},
\hspace{.2in}\mbox{if}\hspace{.1in}\frac{1}{\beta_0} \leq x_1
\leq \frac{1}{\beta}.$$
From this Claim 2 follows.
\vskip 5mm
Now we apply Claim 2 to our situation. There
$$r=\rho-\epsilon,\hspace{.2in}
\beta=\rho-\epsilon-\gamma,$$
where
$$\gamma^2= ||x-x_t||^2 +(\rho-\epsilon)^2
- 2(\rho-\epsilon)||x-x_t||cos\theta.$$

$$\alpha=\gamma-\frac{\rho-\epsilon-s_0}{\gamma}(\gamma^2-||x-x_t||^2
\mbox{sin}^2 \theta)^{\frac{1}{2}}-
\frac{(2(\rho-\epsilon)s_0-s_0^2)^{\frac{1}{2}}||x-x_t||
\mbox{sin}\theta}{\gamma}$$
and
$$\beta_0=\frac{(\rho-\epsilon)^2}{\gamma - \alpha}-\gamma.$$
Then we get, similarily as before
$$\frac{t^{\frac{2}{n+1}}}{n} < x,N(x) > \left(1- ( \frac{\|
x_t\|}{\| x
\|} )^n \right)$$
$$\leq$$
$$\frac{1}{2}
(\frac{|K|}{v_n})^{\frac{2}{n+1}}
(\rho-\epsilon)^{-\frac{n-1}{n+1}}
\frac{1}
{(1+\frac{\| x -x_t
\|^2}{(\rho-\epsilon)^2}-2\frac{\| x -x_t\|\mbox{cos}\theta}
{\rho-\epsilon})^{\frac{1}{2}}(1-\frac{\| x -x_t\|}
{2(\rho-\epsilon)\mbox{cos}\theta})}$$
$$\cdot\{1+\frac{2^{\frac{n+1}{2}}R(\rho-\epsilon)^n}
{v_n}(1-\frac{\gamma}{\rho-\epsilon})^{\frac{n+1}{2}}
(\frac{\gamma}{\rho-\epsilon})^{\frac{n+1}{2}}
\}^{\frac{2}{n+1}},$$
with an $R$ defined accordingly.
\newline
This finishes the proof of Lemma 13 in the case where the
indicatrix is an ellipsoid.
\vskip 5mm
(ii) \underline{Case where the indicatix of Dupin is an elliptic
cylinder}
\newline
Recall that then we have to show that
$$\lim_{t \rightarrow \infty} \frac{1}{n} < x,N(x) >
t^{\frac{2}{n+1}} \left(1- ( \frac{\| x_t \|}{\| x \|} )^n
\right)
=0.$$
We can again assume
(see [S-W]) that the indicatrix is a spherical cylinder i.e. the product of a
$k$-dimensional plane
and a $n-k-1$ dimensional Euclidean sphere of radius $\rho$. Moreover we can
assume that
$\rho$ is arbitrarily large (see also [S-W]).

By Lemma 9 of [S-W] we then have for sufficiently small $s$ and some
$\varepsilon > 0$

$$B_{\rho - \varepsilon} \cap H_s^+
\subseteq K \cap H_s^+.$$
Using similar methods, this implies the claim.

\newpage
\medskip\noindent {\bf Proposition 14}
\par
\begin{it}
Let K be a convex body such that $\partial K$ is $C^3$ and has
strictly positive Gaussian curvature everywhere. Then there is $ \delta_0 > 0$
such that for all $\delta < \delta_0$
$$S(K,\frac{v_{n-1}}{2(n+1)v_n\delta}) \subseteq K_{\delta|K|}.$$
\end{it}

\vskip 3mm
\medskip\noindent {\bf Proof}

As in the proof of Lemma 12, we can choose  $1 > \alpha > 0$ such
that

$$B(0,\alpha) \subseteq K \subseteq
B(0,\frac{1}{\alpha}).$$ Therefore we have for all $x \in \partial K$
$$
1 \geq \hspace{.1in}< \frac{x}{\| x \|} ,N(x) >\hspace{.1in} \geq \alpha^2.$$
Let $R_0 = \mbox{min}_{x \in \partial K,1 \leq i \leq n-1}R_i(x)$, where
$R_i(x)$
is the i-th principal radius of curvature at $x \in \partial K$.
$R_0 > 0$  (see [L2]).
\newline
Let $1 > \varepsilon >0$ be given such that $ \varepsilon  < \mbox{min}
\{\frac{R_0}{2}, 6n\alpha^4\}$  and such that $$(1-\varepsilon)
(1-\frac{\varepsilon}{R_0})^{n-1}(1 -\frac{\varepsilon}{6n\alpha^4})^n
>\frac{1}{2}.$$
By assumption the indicatrix of Dupin exists for all $x \in \partial K$
and is an
ellipsoid. For $x \in \partial K$  given, we can assume that, after an
affine transformation,
the indicatrix at $x$ is a Euclidean sphere.
Let $\sqrt {\rho(x)}$ be the radius of this
Euclidean sphere. Note that for all $x \in \partial K$,
\begin{equation}\rho(x) \geq R_0.
\end{equation}
Then, with the notations used in the proof of Lemma 13,
there exists $s(x) > 0$ such that
$$B_{\rho(x)-\varepsilon} \cap H_{s(x)}^+ \subseteq K \cap H_{s(x)}^+ \subseteq
B_{\rho(x)+\varepsilon} \cap H_{s(x)}^+.$$
Let $s_1=\mbox{min}_{x \in \partial K}s(x)$. $s_1 > 0$ as $\partial K $ is
$C^3$
and compact. Let
\begin{equation}
s_0 = \mbox{min} 
\{s_1, (R_0 -\varepsilon)(1-\frac{1}{1+\frac{2\varepsilon}{3n}})\}.
\end{equation}

Let $\delta_0 >0$ be so small that for all $x \in \partial K$ two
conditions are satisfied;
firstly
\begin{equation}
||x-x_{\delta_0}||<\frac{x}{||x||},N(x)> \hspace{0.1in}\leq
\hspace{0.1in}\frac{s_0}{2},
\end{equation}
where $x_{\delta_0} =[0,x] \cap \partial K_{\delta _0|K|},$
and secondly
$$H_{\delta_0}^+ \cap B_{\rho(x)-\varepsilon} \subseteq
H_{s_0}^+ \cap B_{\rho(x)-\varepsilon},$$
where $H_{\delta_0}$ is the hyperplane through $x_{\delta_0}$
that cuts off exactly $\delta_0|K|$ from $K$.
\vskip 3mm
Suppose now that the above Proposition is not true. Then there is $\delta
<\delta_0$
and $x_s \in \partial S(K,\frac{v_{n-1}}{2(n+1)v_n \delta})$ such that $x_s
\notin
K_{\delta|K|}$. Let $x \in \partial K$ be such
that $x_s \in[0,x]$. We also can assume
that the  indicatrix of Dupin at x is a Euclidean ball with radius
$\sqrt{\rho(x)}$.
We choose $x_{\delta} \in \partial K_{\delta|K|}$ such that $x_{\delta}
\in[0,x]$. Then $$||x-x_{\delta}||<\frac{x}{||x||},N(x)>
\hspace{0.1in}>\hspace{0.1in}
||x-x_{s}||<\frac{x}{||x||},N(x)>.$$
By construction
$$\delta|K| = |K \cap H_{\delta}^+| \geq |B_{\rho(x)  -\varepsilon}
\cap H_{\delta}^+|$$
$$\geq \mbox{min}_{H \in {\cal H}} |B_{\rho(x)  -\varepsilon} \cap H^+|=
|B_{\rho(x)  -\varepsilon} \cap H_0^+|,$$
where $$ {\cal H} = \{H: H \hspace{0.03in}\mbox{is hyperplane through}
\hspace{0.03in} x_\delta,\hspace{0.03in}
x \in H^+\}.$$
For the height $h$ of this cap $|B_{\rho(x)  -\varepsilon} \cap H_0^+|$
of $B_{\rho(x) -\varepsilon}$ of minimal volume we have
$$h \geq ||x-x_{\delta}||<\frac{x}{||x||},N(x)>
(1-\frac{||x-x_{\delta}||}{2<\frac{x}{||x||},N(x)>(\rho(x) -\varepsilon)}).$$
Using (10), (11) and (12)  we get
$$h \geq ||x-x_{\delta}||<\frac{x}{||x||},N(x)>
(1-\frac{\varepsilon}{6n\alpha^4}).$$

The volume of a cap of a Euclidean ball with radius r and height h can be
estimated from below by
$$ \geq \frac{v_{n-1} 2^{\frac{n+1}{2}}r^{\frac{n-1}{2}}
h^{\frac{n+1}{2}}}{n+1}(1-\frac{h}{2r})^{\frac{n+1}{2}}.$$
Therefore
$$\delta|K| \geq$$
$$\frac{v_{n-1}}{n+1} 2^{\frac{n+1}{2}}\rho(x)^{\frac{n-1}{2}}
(||x-x_{\delta}||<\frac{x}{||x||},N(x)>)^{\frac{n+1}{2}}
(1-\frac{\varepsilon}{R_0})
^{\frac{n^-1}{2}}
(1-\frac{\varepsilon}{6n\alpha ^4})^{n}
,$$
where we have used again (10), (11), (12) and the fact that
$<\frac{x}{||x||},N(x)> \geq \alpha ^2$.
Thus
\begin{equation}
(||x-x_{s}||<\frac{x}{||x||},N(x)>)^{\frac{n+1}{2}} <
\frac {(n+1)\delta |K|\rho (x)^{-\frac{n-1}{2}}}
{v_{n-1} 2^{\frac{n+1}{2}}
(1-\frac{\varepsilon}{6n\alpha ^4})^n
(1-\frac{\varepsilon}{R_0})^{\frac{n-1}{2}}}.
\end{equation}
As, with the notations of Lemma 13 and using (10), (11) and (12)
$$\frac{c}{\rho(x) +\varepsilon} \geq 1-
\frac{||x-x_{s}||<\frac{x}{||x||},N(x)>}{\rho(x) +\varepsilon} >
\frac{1}{1+\frac{2\varepsilon}{3n}},$$
the estimate from below from Lemma 13 for $\frac{|K||K^x|}{v_n^2}$ holds
for $x=x_s$
and we get (see p. 23)
$$\frac{|K||K^{x_s}|}{v_n^2} > \frac{(1-\varepsilon)|K| \rho(x)^{-n}
(1+\frac{\varepsilon}{ \rho(x)})^{-n}}{2^{\frac{n+1}{2}}v_n
(\frac{c}{\rho(x) +\varepsilon})^{\frac{n+1}{2}}
(1-\frac{c}{\rho(x) +\varepsilon})^{\frac{n+1}{2}}}.$$
Now notice that
$$\frac{c}{\rho(x) +\varepsilon} \leq 1$$
and
$$1-\frac{c}{\rho(x) +\varepsilon} \leq
\frac{||x-x_{s}||<\frac{x}{||x||},N(x)>}{\rho(x)
(1+\frac{\varepsilon}{\rho(x)}}.$$
Therefore (13) implies that
$$\frac{|K||K^{x_s}|}{v_n^2}> (1-\varepsilon)
(1-\frac{\varepsilon}{R_0})^{{n-1}}(1 -\frac{\varepsilon}{6n\alpha^4})^n
\quad \frac{v_{n-1}}{(n+1)\delta v_n}.$$
This is a contradiction.

\vskip 60mm
\medskip\noindent Mathieu Meyer
\newline
Universit\'{e} de Marne-la-Vall\'{e}e
\newline
Equipe d'Analyse et de Math\'{e}matiques Appliqu\'{e}es,
\newline
5 rue de la Butte verte
\newline
93166 Noisy-le-Grand Cedex, France
\newline
e-mail: meyer@math.univ-mlv.fr
\vskip 3mm
\medskip\noindent Elisabeth Werner
\newline
Department of Mathematics
\newline
Case Western Reserve University
\newline
Cleveland, Ohio 44106, U.S.A.
\newline
e-mail: emw2@po.cwru.edu
\newline
and
\newline
Universit\'{e} de Lille 1
\newline
UFR de Math\'{e}matiques
\newline
59655 Villeneuve d'Ascq, France

\end{document}